\documentclass[12pt, a4paper]{amsart}
\usepackage{latexsym,amssymb,amsmath}
\newtheorem{thm}{Theorem}[section]

\newtheorem{lem}[thm]{Lemma}
\newtheorem{prop}[thm]{Proposition}
\newtheorem{example}[thm]{Example}

\begin{document}
\title{Curves of genus 3 over small finite fields}
\author{Jaap Top}
\address{IWI, Rijksuniversiteit Groningen, 
Postbus 800, NL-9700 AV Groningen, The Netherlands}
\email{top@math.rug.nl}
\maketitle

\section{Introduction}
The maximal number of rational points that a
(smooth, geometrically irreducible) curve of genus $g$ over
a finite field ${\mathbb F}_q$ of cardinality $q$ can have, is denoted by $N_q(g)$.
The interest in this number, particularly for fixed $q$ as a function in $g$,
arose primarily during the last two decades from applications to error
correcting codes \cite{L}, \cite{TGZ}, \cite{LG}. A lot of results on
$N_q(g)$ for fixed $q$ and `$g$ large' are discussed in the second part of
J-P.~Serre's 1985 Harvard lectures \cite{Se85}.

The first part of these Harvard lectures reverses the roles of $q$ and $g$:
fix the genus $g$ and study the resulting function $N_q(g)$ of $q$.
In terms of coding theory, this implies one does not consider
asymptotic results (in terms of increasing lengths of the codes), but
rather one puts a constraint on the `complexity' of the curves used
for constructing codes.

It is a classical result that $N_q(0)=q+1$: any curve of genus $0$
containing a rational point, is isomorphic to the projective line.
The determination of $N_q(1)$ is the work of Max Deuring, published
in 1941 (see \cite{Deuring} and also \cite{Wat}). Next, $N_q(2)$ is
computed by
J-P.~Serre (1983) \cite{Se1}, \cite{Se2}, \cite{Se85}; see also
\cite{Shabat}.

For $g\geq 3$, no general formula for $N_q(g)$ seems to be known. The tables
\cite{G-V} describe what is known about $N_q(g)$ for $g\leq 50$ and
$q\in\{2,3,4,8,9,16,27,32,64,81,128\}$. J-P.~Serre in \cite[\S~4]{Se2}
and \cite[p.~64-65]{Se85} presents the values of $N_q(3)$ for
all $q\leq 25$. This table is slightly extended in \cite[p.~164]{ATants}.

The goal of the present note is twofold. On the one hand, we prove
a `guess' expressed by J-P.~Serre on page~66 of \cite{Se85}: 
\begin{prop}\label{F-special}
If $C/{\mathbb F}_q$ is a curve of genus $3$ with the property
$\#C({\mathbb F}_q)>2q+6$, then $q\in\{8,9\}$ and $C$ is isomorphic over 
${\mathbb F}_q$ to one of the following two curves.
\begin{enumerate}
\item ($q=8$). The plane curve over ${\mathbb F}_8$ given by $x^4+y^4+z^4+x^2y^2+y^2z^2
+x^2z^2+x^2yz+xy^2z+xyz^2=0$, which has exactly $24$ rational points over
${\mathbb F}_8$;
\item ($q=9$). The quartic Fermat curve over ${\mathbb F}_9$ given by 
$x^4+y^4+z^4=0$, which has exactly $28$ rational points over
${\mathbb F}_9$.
\end{enumerate}
\end{prop}

The proof uses work of K.-O.~St\"{o}hr and J.F.~Voloch \cite{S-V}, which allow one to
translate the asumption into a property of the endomorphism ring of
the jacobian of the curve. Also, a result of A. Hefez and J.F.~Voloch \cite{H-V}
classifying plane quartics with the property that every tangent line is in
fact a flex line, is used. A large part of the proof can be traced back
to \cite{Se85}.

Secondly, some results on determining $N_q(3)$ for relatively small $q$ are
discussed. For this, various known upper bounds are recalled, and some
explicit constructions of genus $3$ curves on which it is easy to
determine the number of rational points are presented. Using this,
the following extension to Serre's table was found.
\begin{prop}\label{table}
The maximum number $N_q(3)$ of rational points on a curve of genus $3$
over a finite field with $q$ elements is, for $q<100$, given in the following
table.
$$
\begin{array}{c||c|c|c|c|c|c|c}
q      & 2 & 3 & 4 & 5 & 7 & 8 & 9   \\ 
N_q(3) & 7 & 10& 14& 16& 20& 24& 28  \\
\hline
q & 11 & 13 & 16 & 17 & 19 & 23 & 25  \\
N_q(3)  & 28 & 32 & 38 & 40 & 44 & 48 & 56  \\
\hline
q     & 27 & 29 & 31 &32 & 37 & 41 & 43  \\
N_q(3) & 56 & 60&  62&64 & 72 & 78&80\\
\hline
q     & 47 & 49 & 53 &59 & 61 & 64 & 67 \\
N_q(3) &87&92&96&102&107&
113&116\\
\hline
q     &71 & 73 & 79 & 81 & 83 & 89 & 97\\
N_q(3) &120&122&131&136&136&144&155
\end{array}
$$ 
\end{prop}

Proposition~\ref{F-special} is proven in Section~\ref{Frob-s} below,
and Proposition~\ref{table} in Sections~\ref{UB} and \ref{Curves}.

This note originated in three lectures on the results
of \cite{ATants} I gave in the
Mathematics Departments of Utrecht, Eindhoven and
Nijmegen between October 2002 and January 2003. I thank the
colloquium organizers Gunther Cornelissen and Johan van de Leur
in Utrecht, the organizer Arjeh Cohen of the EIDMA seminar
{\sl Combinatorial Theory}, and Jozef Steenbrink in Nijmegen for
offering the opportunity to speak about these topics.

\section{Upper bounds}\label{UB}
Let $q$ be the cardinality of some finite field ${\mathbb F}_q$.
Suppose $C/{\mathbb F}_q$ is a (smooth, complete, absolutely irreducible)
curve of genus $3$. We write $N$ for the cardinality of the set
$C({\mathbb F}_q)$ of ${\mathbb F}_q$-rational points on $C$.
Finally, we put $m:=[2\sqrt{q}]$, the largest integer which is
$\leq 2\sqrt{q}$. The following bounds for the number $N$ are known.
\begin{prop}\label{ubounds}
With notations as above, one has
\begin{itemize}
\item[(a)] $N\leq q^2+q+1$;
\item[(b)] $N\leq 2q+6$ except for $q=8$ and for $q=9$;
\item[(c)] $N\leq q+1+3m$;
\item[(d)] $N\leq q+3m-1$ if $q=a^2+1$ for some integer $a$;
\item[(e)] $N\leq q+3m-1$ if $q=a^2+2$ for some integer $a\geq 2$;
\item[(f)] $N\leq q+3m-2$ if $q=a^2+a+1$ for some integer $a$;
\item[(g)] $N\leq q+3m-2$ if $q=a^2+a+3$ for some integer $a\geq 3$.
\end{itemize}
\end{prop}

The first assertion here follows from the well known fact that
there are only two possibilities for a curve of genus $3$. 
It can be hyperelliptic, which means that
its canonical morphism to ${\mathbb P}^2$ has degree $2$ and
has as image a rational curve. This implies that $N\leq 2q+2$,
which is less than $q^2+q+1$ since $q\geq 2$. If the
curve is not hyperelliptic, then the canonical morphism defines
an isomorphism to a plane quartic curve. In this case of course
$N$ is at most the cardinality of ${\mathbb P}^2({\mathbb F}_q)$,
which is $q^2+q+1$.
  
\begin{example}\label{example} \rm
The equation $y^2z^2+yz^3+xy^3+x^2y^2+x^3z+xz^3=0$ defines a smooth quartic
over ${\mathbb F}_2$ containing all $7$ points of ${\mathbb P}^2({\mathbb F}_2)$.
Hence $N_2(3)=7$, as is asserted in various tables such as
\cite{Se1}, \cite[p.~64-65]{Se85}, \cite[p.~71]{LG}, \cite{G-V}.
Note that in \cite[p.~64-65]{Se85}, Serre remarks that ``{\sl this list is not
entirely guaranteed $\ldots$}''. Indeed, the curve
$x^4+y^4+z^4+x^2y^2+y^2z^2+x^2z^2+x^2yz+xy^2z=0$ over ${\mathbb F}_2$ given
there, contains $0$ rather than $7$ rational points over ${\mathbb F}_2$.
\end{example}

Item (b) in Proposition~\ref{ubounds} follows from a beautiful geometric
argument presented in \cite{S-V} by K.-O.~St\"{o}hr and J.F.~Voloch;
see also \cite[p.~64-66]{Se85}.
We recall it here, since it is relevant for the proof of
Proposition~\ref{F-special} given in Section~\ref{Frob-s} below.
The assertion is true for hyperelliptic curves, hence we can
and will assume that $C$ is a plane quartic given by $F=0$.
Writing $F_x$, et cetera, for the partial derivative 
$\frac{\partial F}{\partial x}$, the line $T_p$ in ${\mathbb P}^2$ tangent
to $C$ in a point $p\in C$ is given by
$xF_x(p)+yF_y(p)+zF_z(p)=0$. Write $\pi:C\to C$ for the Frobenius
morphism over ${\mathbb F}_q$, which raises all coordinates to the
$q$th power. Obviously, one has for a point $p\in C$ that
$$p\in C({\mathbb F}_q)
\Leftrightarrow
\pi(p)=p
\Rightarrow
\pi(p)\in T_p.
$$
The latter condition can be written as
$$
x(p)^qF_x(p)+y(p)^qF_y(p)+z(p)^qF_z(p)=0
$$
where $x(p), y(p), z(p)$ are the projective coordinates of the point $p$.
Hence if we denote by $C_q$ the plane curve of degree $q+3$ given by the equation
$x^qF_x+y^qF_y+z^qF_z=0$, then it follows that
$$C({\mathbb F}_q)\subseteq C\cap C_q.$$
Now two possibilities arise.

If $C\subset C_q$, then $C/{\mathbb F}_q$ is called
{\sl Frobenius non-classical}. It means that {\em every} $p\in C$ has the
property that $\pi(p)\in T_p$. Proposition~\ref{F-special} asserts
in particular, that this happens only for $q=8$ and for $q=9$,
which we have excluded here.

If $C\not\subset C_q$, then the number of intersection points of $C$
and $C_q$, counted with multiplicities, equals $4\cdot(q+3)$ by
B\'{e}zout's theorem \cite[Ch.~I,~Cor.~7.8]{HAG}. Note that at every point $p\in C({\mathbb F}_q)$,
the line $T_p$ intersects both $C$ and $C_q$ at $p$ with
multiplicity $\geq 2$. Hence these points also have
intersection multiplicity $\geq 2$ in $C\cap C_q$. It
follows that $2N\leq 4\cdot(q+3)$, which is what we wanted to prove. 

Item (c) in Proposition~\ref{ubounds} is J-P.~Serre's
refinement \cite{Se1} of the classical Hasse-Weil bound.

Items (d), (e), (f) and (g) can be found in Lauter's paper
\cite[Thm.~1 and 2]{LauS}. Basically, the proof extends an idea
exposed in \cite[p.~13-15]{Se85}, which Serre attributes to
A.~Beauville. 

\section{Frobenius non-classical quartics}\label{Frob-s}
In this section we present a proof of Proposition~\ref{F-special}.
So we suppose $C/{\mathbb F}_q$ is a curve of genus $3$ satisfying
$\#C({\mathbb F}_q)>2q+6$. Then in particular $C$ is not hyperelliptic,
hence we may assume that $C$ is a plane quartic curve. The St\"{o}hr-Voloch
result as given in the previous section shows that $C$ is 
Frobenius non-classical. The tangent line $T_x$ at a general point
$x\in C$ then intersects $C$ in $x$ (with multiplicity $\geq 2$),
also in $\pi(x)\neq x$, and hence in a unique fourth point $\phi(x)$.
The map $\phi:\;C\to C$ defines a morphism. It is non-constant because
a smooth curve of positive genus cannot have the property that all
its tangent lines meet in a single point (see \cite[Ch.~IV,~Thm.~3.9]{HAG}).

If $x,y$ are two points on $C$ with tangent lines $T_x,T_y$ given by
$\ell=0$ and $m=0$ respectively, then $f:=\ell/m$ defines a function
on $C$ with divisor $\mbox{div}(f)=2(x-y)+\pi(x-y)+\phi(x-y)$.
In particular, this implies that in the endomorphism ring of the
Jacobian $J(C)$ of $C$, the relation $2+\pi+\phi=0$ holds. We will
distinguish two cases.

(1). If $\phi\;:C\to C$ is separable, then the Hurwitz formula implies
that $\phi$ is an isomorphism. As a consequence, the eigenvalues of
$\phi$ acting on a Tate module of $J(C)$ are roots of unity. Let
$\zeta$ be such an eigenvalue. The equality $2+\pi+\phi=0$ shows
that $-2-\zeta$ is an eigenvalue of Frobenius, hence
$(-2-\zeta)(-2-\zeta^{-1})=q$. It follows that $q=|2+\zeta|^2\leq 9$.
Moreover, since $\zeta+\zeta^{-1}=(q-5)/2$ is integral, it follows that $q$ is odd. 
So $q\in\{3,5,7,9\}$.
If $q=3$ then $\zeta$ would be a primitive third root of unity $\omega$
and the eigenvalues of $\pi$ would be $-2-\omega$ and $-2-\omega^{-1}$,
each with multiplicity $3$. Hence the trace of Frobenius would be $-9$
which implies $\#C({\mathbb F}_3)=13$. However, it is well known that
over ${\mathbb F}_3$, a genus 3 curve cannot have more than $10$ rational
points (\cite[\S~4]{Se2}, \cite{G-V} and also Proposition~\ref{ubounds}(f) above).
If $q=5$ or $q=7$, then every eigenvalue of Frobenius is $-2-\zeta$ for
a primitive fourth respectively sixth root of unity $\zeta$. This implies that the
only eigenvalue of $\phi^{(q-1)/2}$ is $-1$, hence also $\phi^{(q-1)/2}$
acting on the space of regular $1$-forms on $C$ has as only eigenvalue $-1$.
Therefore $C$ is hyperelliptic, contrary to our assumption.
The remaining case is that $q=9$ and $\zeta=1$. It means that $\phi$
is the identity map, hence every line tangent to $C$ is in fact a flex line.
As is shown in the paper by Hefez and Voloch \cite{H-V}, this implies
in our situation that $C$ is isomorphic to the Fermat curve given by
$x^4+y^4+z^4=0$. We also have that every eigenvalue of Frobenius
equals $-3$, hence $\#C({\mathbb F}_9)=9+1+6\cdot3=28$. 

(2). The remaining case is that $\phi\;:C\to C$ is inseparable.
Then both $\pi$ and $\phi$ are inseparable and hence so is $2=-\pi-\phi$.
This implies that $q$ is even. Write $\pi_2$ for the map
that raises coordinates to their second power. Since $\phi$ is inseparable,
one can write $\phi=\psi\pi_2^n$ for some $n>0$ and $\psi$ separable.
Put $\pi=\pi_2^m$, then with $k:=\mbox{min}(m,n)$ one finds
$2+(\pi_2^{(m-k)}+\psi\pi_2^{(n-k)})\pi_2^k=0$. Comparing degrees (of
isogenies between Jacobians) it follows that $2^{3k}\leq 2^6$ hence $k=1$
or $k=2$. Therefore, considering a pair $\lambda,\mu$ of eigenvalues
of $\pi$ and $\phi$, with $2+\lambda+\mu=0$ and $|\lambda|^2=2^{k}\leq |\mu|^2$,
one finds, using that $|\mu|^2=(2+\lambda)(2+\bar{\lambda})$ is a power of $2$,
that $q\leq 16$. On the other hand, a result of \cite{H-V} says that
any Frobenius non-classical curve of degree $d$ over ${\mathbb F}_q$ contains
exactly $d(q-d+2)$ rational points. In our case, this equals $4q-8$, which
is $>2q+6$ only when $q>7$. On the other hand, this number cannot exceed
the Hasse-Weil bound $q+1+6\sqrt{q}$, hence we must have $q\leq 9$.

It follows that $q=8$. In this case, the possibilities for the pairs $\lambda,\mu$
show that either the eigenvalues of Frobenius are $-2\pm2i$ (each with multiplicity $3$,
which would yield the number of points $8+1+12=21<2q+6$,
or the eigenvalues $\lambda=(-5\pm \sqrt{-7})/2$ (each with multiplicity $3$),
which yield as the number of points $8+1+15=24=4q-8$. Such a $\lambda$ corresponds
to $\mu=(1\mp \sqrt{-7})/2$. Note that $\mu^3=\lambda$. Hence on an abelian
$3$-fold $A$ over ${\mathbb F}_2$ where the Frobenius $\pi_2$ satisfies
$\pi_2^2-\pi_2+2=0$, we will have $\pi_2^3+\pi+2=0$. Over ${\mathbb F}_8$ this
$A$ will be isogenous to our $J(C)$. It follows that $C$ can be defined
over ${\mathbb F}_2$, and in fact is given by a plane quartic curve there
with the property that for every point $x$ on it, $x,\pi_2(x)$ and $\pi_2^3(x)$
are collinear. In terms of the coordinates $(\xi,\eta,\mu)$ of such a point,
it means that the three vectors $(\xi,\eta,\mu),\;(\xi^2,\eta^2,\mu^2)$
and $(\xi^8,\eta^8,\mu^8)$ are linearly dependent. Thus, $(\xi,\eta,\mu)$
satisfies the degree $11$ homogeneous polynomial
$$\det\left(
\begin{array}{ccc}\xi&\xi^2&\xi^8\\
\eta&\eta^2&\eta^8\\
\mu&\mu^2&\mu^8\end{array}\right)=0.$$
Clearly, this polynomial is divisible by all of the $7$ homogeneous degree $1$
polynomials over ${\mathbb F}_2$. So the curve(s) we look for will be given
by the remaining factor of degree $11-7=4$, which is
$$\xi^4+\eta^4+\mu^4+\xi^2\eta^2+\eta^2\mu^2+\xi^2\mu^2
+\xi^2\eta\mu++\xi\eta^2\mu+\xi\eta\mu^2=0.$$
This indeed defines a smooth quartic, which by construction is
Frobenius non-classical over ${\mathbb F}_8$.

This finishes the proof of Proposition~\ref{F-special}.

\section{Constructions of curves with many points}\label{Curves}
In this section we provide the information needed to verify Proposition~\ref{table}.
First of all, in most cases Proposition~\ref{ubounds} shows that the values presented in our table
are upper bounds for $N_q(3)$. To be precise, we mention in the following table
which of the items (a)--(f) in Proposition~\ref{ubounds} works for a given $q$.

\medskip
\begin{tabular}{|r|ccccccc|}\hline
$q$:& 2&3&4&5&7&8&9\\
item(s):&(a,d)&(f)&(b)&(b,d)&(b,f)&(c)&(c)\\\hline
$q$:&11&13&16&17&19&23&25\\
item(s):&(b,e)&(b,f)&(b)&(b,d)&(b,c)&(g)&(b,c)\\\hline
$q$:&27&29&31&32&37&41&43\\
item(s):&(e)&(c)&(f)&--&(d)&(c)&(f)\\\hline
$q$:&47&49&53&59&61&64&67\\
item(s):&(c)&(c)&(c)&(g)&(c)&(c)&(c)\\\hline
$q$:&71&73&79&81&83&89&97\\
item(s):&(c)&(f)&(c)&(c)&(e)&(c)&(c)\\\hline
\end{tabular}

\medskip
For $q=32$, Proposition~\ref{ubounds} only gives the bound $N\leq 66=q+1+3m$.
However, the maximum in this case is $64$, as is asserted in the tables \cite{G-V}.
A curve over ${\mathbb F}_{32}$ containing $64$ rational points is given
in the thesis \cite{Sem}. To show that this is indeed the maximum, one may use the following fact
due to J-P.~Serre.
\begin{lem}[Serre]
Let $q>1$ be a power of a prime number and $m:=[2\sqrt{q}]$. If $4q-m^2\leq 11$
then $N_q(3)\leq q+3m-1$.
\end{lem}
Proof: The fact that under the given conditions no genus $3$ curve with exactly
$q+1+3m$ points exists, follows from \cite[App.~7.1]{LauS}. The fact that
no curve of genus $3$ with eactly $q+3m$ points exists, is proven in \cite{Se2}.
This implies the lemma. \hfill{$\Box$}

\vspace{\baselineskip}
What remains, is to show that for each $q<100$ the given upper bound is sharp, i.e.,
there exists a genus $3$ curve with that number of rational points. For each $q\leq 25$,
this is done by Serre \cite[p.~64-65]{Se85} (note Example~\ref{example} given above,
which seems to be the only `misprint' in Serre's table).
We will now discuss the cases $25<q<100$.

In the case $q=27$, a genus $3$ curve with exactly $56$ points was constructed by
Van der Geer and Van der Vlugt \cite{G-V}.

For $q\in\{29,49,53,67,71,89\}$, a curve with $q+1+3m$ points exists. In fact,
such a curve is found in the family
$$C_{\lambda}:\;x^4+y^4+z^4=(\lambda+1)(x^2y^2+y^2z^2+x^2z^2)$$
as follows from the table given in \cite[\S~4.2]{ATants}. 
Specifically, one finds $\#C_2({\mathbb F}_{29})=60$ and
$\#C_{-1}({\mathbb F}_{49})=92$ and
$\#C_2({\mathbb F}_{53})=96$ and
$\#C_{30}({\mathbb F}_{67})=116$
and $\#C_{37}({\mathbb F}_{71})=120$ and
$\#C_{13}({\mathbb F}_{89})=144$.
The same
family $C_{\lambda}$ yields an optimal curve for $q=43$, with $q+1+3(m-1)=80$ rational points,
namely $C_{10}$.

The prime numbers $q\in\{31,61,73,79,97\}$ can be treated using the family
$$D_{a,b}:\; x^3z+y^3z+x^2y^2+axyz^2+bz^4=0$$
of curves with a noncyclic automorphism group of order $6$.
In this family one finds, for instance, $\#D_{4,2}({\mathbb F}_{31})=62$ and
$\#D_{29,34}({\mathbb F}_{61})=107$ and
$\#D_{2,48}({\mathbb F}_{73})=122$ and
$\#D_{11,8}({\mathbb F}_{79})=131$ and
$\#D_{56,79}({\mathbb F}_{97})=155$.

The case $q=32$ we already discussed.

For $q=37$ and $q=83$ one may proceed as follows.
Consider the elliptic curves
$$E_{\lambda}:\;y^2=x(x-1)(x-\lambda)$$
and
$$E'_{a,b}:\; y^2=x(x^2+ax+b).$$
Denote by $T$ resp. $T'$ the point $(0,0)$ on $E_\lambda$ and $E'_{a,b}$, respectively.
Let $P=(1,0)\in E_{\lambda}$ and let $Q\neq T'$ be one of
the points with $y$-coordinate $0$ on $E'_{a,b}$.
The triples $(0,T',T')$ and $(T,T',0)$ and $(P,Q,Q)$ generate a
(rational) subgroup $H$ of order $8$ in $E_\lambda\times E'_{a,b}\times E'_{a,b}$. 
The quotient $(E_\lambda\times E'_{a,b}\times E'_{a,b})/H$ is an
abelian variety. From \cite[Prop.~15]{HLP} it follows that this is in
fact the jacobian of a genus $3$ curve $C$ over ${\mathbb F}_q$,
provided that $\lambda(\lambda-1)(a^2\lambda-4b)$ is a nonzero square
in ${\mathbb F}_q$. Under this condition, one has
$$\#C({\mathbb F}_q)=\#E_{\lambda}({\mathbb F}_q)+2\#E'_{a,b}({\mathbb F}_q)-2q-2.$$
In fact, the curve $C$ here can be given by the equation
$$(\lambda-1)\left(\lambda x^4+by^4+bz^4+\lambda a x^2y^2+\lambda a x^2z^2\right)
=(\lambda a^2-2b(\lambda+1))y^2z^2.$$
The choice $\lambda=7$, $a=0$ and $b=2$ turns out to give
a curve $C$ over ${\mathbb F}_{37}$ with exactly $72$ rational points. 
Similarly, $\lambda=5$ and $(a,b)=(4,2)$ yields a genus $3$ curve over
${\mathbb F}_{83}$ containing $136$ rational points.

For $q=41$, consider the family
$$X_{a,b}:\;x^2y^2+y^2z^2+x^2z^2+a(x^3y+y^3z+xz^3)+b(x^3z+xy^3+yz^3)=0$$
of curves with an automorphism of order $3$. Here one finds the example
$\#X_{-7,8}({\mathbb F}_{41})=78$.

In case $q=59$ one can use the family
$$Y_{a,b}:\;(3x^2+y^2)^2+ax(x^2-y^2)z+bz^4=0,$$
which are curves with a noncyclic group of automorphisms (of order $6$)
generated by $(x,y,z)\mapsto (x,-y,z)$ and $(x,y,z)\mapsto(-x-y,-3x+y,2z)$.
One finds for example $\#Y_{4,6}({\mathbb F}_{59})=102$.

For $q=64$ we quote \cite{G-V}, which states that a genus $3$ curve
over ${\mathbb F}_{64}$ containing $113$ rational points is
given in \cite{Wi}. Similarly, for the case $q=81$ the
reference \cite{Wi} provides, according to \cite{G-V}, an example with 
$136$ rational points.

It should be remarked here that the families $C_{\lambda}$, $D_{a,b}$, $X_{a,b}$, and $Y_{a,b}$
came up in a rather natural way when searching for curves:
namely, in many cases we looked for a maximal curve, which means
one with $q+1+3m$ rational points. In such a case, if $E$ is an
elliptic curve over ${\mathbb F}_q$ containing $q+1+m$ rational points,
the Jacobian of the curve we try to find is isogenous to
$E\times E\times E$. As is explained in the Appendix of \cite{LauS}, the principal 
polarization on this Jacobian
can be interpreted as a rank $3$ unimodular indecomposable hermitian module
over  ${\mathbb Z}[(-m+\sqrt{m^2-4q})/2]$. The possible group(s) of isometries
of such modules yield the possible automorphism groups of the curves we try to find.
It turned out that in most cases this gave sufficiently many restrictions
on the curves to be actually able to find them.

\end{document}